\documentstyle[11pt,twoside]{article}
\title{\bf A $(p,\nu)$-extension of the Appell function $F_1(\cdot)$ and its properties}
\author{\sc S.A. Dar$^a$ and R.B. Paris$^b$\\
\\
${}^a\!$ {\em Department of Applied Sciences and Humanities, Faculty of Engineering }\\
{\em and Technology, Jamia Millia Islamia, New Delhi, 110025, India}\\ 
{\em E-Mail: showkatjmi34@gmail.com}\\
${}^b\!$ {\em Division of Computing and Mathematics, Abertay University,}\\
{\em Dundee DD1 1HG, UK}\\
{\em E-Mail: r.paris@abertay.ac.uk}
}
\begin{document}
\newcommand{\bee}{\begin{equation}}
\newcommand{\ee}{\end{equation}}
\def\f#1#2{\mbox{${\textstyle \frac{#1}{#2}}$}}
\def\dfrac#1#2{\displaystyle{\frac{#1}{#2}}}
\newcommand{\fr}{\frac{1}{2}}
\newcommand{\fs}{\f{1}{2}}
\date{}
\maketitle
\pagestyle{myheadings}
\markboth{\hfill \it S.A. Dar and R.B. Paris  \hfill}
{\hfill \it A $(p,\nu)$ extension of the Appell function \hfill}
\begin{abstract}
In this paper, we obtain a $(p,v)$-extension of the Appell hypergeometric function $ F_{1}(\cdot)$,
together with its integral representation, by using the extended Beta function $B_{p,v}(x,y)$ introduced in \cite{AB30}.  Also, we give some of its main properties, namely the Mellin transform, a differential formula, recursion formulas and a bounded inequality. In addition, some new integral representations of the extended Appell  function$ F_{1,p,v}(\cdot)$ involving Meijer's $G$-function are obtained.
\vspace{0.4cm}

\vspace{0.3cm}

\end{abstract}

\vspace{0.3cm}

\noindent $\,$\hrulefill $\,$

\vspace{0.2cm}

\begin{center}
{\bf 1. \  Introduction and Preliminaries}
\end{center}
\setcounter{section}{1}
\setcounter{equation}{0}
\renewcommand{\theequation}{\arabic{section}.\arabic{equation}}
In the present paper, we employ the following notations:  
\[{\bf N}:=\{1,2,...\},~~{\bf N}_{0}:={\bf N}\cup\{0\},~~{\bf Z}_{0}^{-}:={\bf Z}^{-}\cup\{0\},\]
where the symbols ${\bf N}$ and ${\bf Z}$ denote the set of integer and natural numbers. Also, for the sake of conciseness, we use the following notation \cite[p.33]{CR}
\[{\bf C}_>:=\{p\in {\bf C} : \Re(p)>0\},\]
where, as usual, the symbol ${\bf C}$  denotes the set of complex numbers

In the available literature on hypergeometric series, the hypergeometric series and its generalizations appear in various branches of mathematics associated with applications. This type of series appears very naturally in quantum field theory, in particular
in the computation of analytic expressions for Feynman integrals. Such integrals can be obtained and computed in different ways which may lead to identities for Appell series; see \cite{AB333}. On the other hand, the application of known relations for Appell series may lead to simplifications, help to solve problems or lead to greater insight in quantum field theory. 

The Appell hypergeometric series are a natural two-variable extension of hypergeometric series, which are treated in detail in ´ Erd\'elyi {\it et al.} \cite{AB17}. In this paper, we highlight some of the most important properties and relations satisfied by these series. In the following, we follow to a great extent the expositions from the classical
texts of Bailey \cite{AB7}, and Slater \cite{AB41} (both contain a great amount of material on hypergeometric series). There are four types of Appell functions denoted by
$F_{1}, F_{2}, F_{3}, F_{4}$; in the present study we shall only be concerned with the first Appell function $F_{1}$  given by \cite[(16.13.1)]{AB26}
\begin{eqnarray}
F_{1}(b_{1},b_{2},b_{3};c_{1};x,y)&=&\sum_{m,n=0}^{\infty}\frac{(b_{1})_{m+n}(b_{2})_{m}(b_{3})_{n}}{(c_{1})_{m+n}}\,
\frac{x^{m}}{m!}\frac{y^{n}}{n!}\nonumber\\
&=&\sum_{n,m=0}^{\infty}\frac{(b_{2})_{m}(b_{3})_{n}~B(b_{1}+m+n,c_{1}-b_{1})}{B(b_{1},c_{1}-b_{1})}\,\frac{x^{m}}{m!}\frac{y^{n}}{n!},\label{ES9}
\end{eqnarray}
%where $(\Re(\gamma)>\Re(\beta)>0,~|z|<1).$\\
where $|x|<1$, $|y|<1$. Here $(\lambda)_{\upsilon}~(\lambda, \upsilon\in {\bf C})$ denotes the Pochhammer symbol (or the shifted factorial,
 since $(1)_{n}=n!)$  defined by
\[
(\lambda)_{\upsilon}:=\frac{\Gamma(\lambda+\upsilon)}{\Gamma(\lambda)}=\left\{\begin{array}{ll}1, & (\upsilon=0,\ \lambda\in{\bf C}\backslash\{0\})\\
\lambda(\lambda+1)...(\lambda+n-1),  & (\upsilon=n\in{\bf N},~\lambda\in{\bf C}) 
 \end{array}\right.
\]
and $B(\alpha, \beta)$ denotes the classical Beta function  defined by \cite[(5.12.1)]{AB26}
\begin{equation}\label{ES5}
B(\alpha, \beta)=\left\{\begin{array}{ll} \displaystyle{\int_{0}^{1}t^{\alpha-1}(1-t)^{\beta-1}dt}, & (\Re(\alpha)>0,  \Re(\beta)>0)\\
\\
\displaystyle{\frac{\Gamma(\alpha)\Gamma(\beta)}{\Gamma(\alpha+\beta)}}, & ( (\alpha,\beta)\notin {\bf Z}_{0}^{-}). \end{array}\right.
\end{equation}
An integral representation of $F_{1}(\cdot)$ is given by \cite[(16.15.1)]{AB26}
\begin{equation}\label{ES10}
F_{1}(b_{1},b_{2},b_{3};c_{1};x,y)=\frac{\Gamma(c_{1})}{\Gamma(b_{1})\Gamma(c_{1}-b_{1})}\int_{0}^{1} t^{b_{1}-1}(1-t)^{c_{1}-b_{1}-1}(1-xt)^{-b_{2}}(1-yt)^{-b_{3}}dt,
\end{equation}
where $\Re(c_{1})>\Re(b_{1})>0$, $|\arg(1-x)|<\pi$ and $|\arg(1-y)|<\pi$.
 
A natural generalization of the Gauss hypergeometric series ${}_{2}F_{1}$ is the general hypergeometric series ${}_{p}F_{q}$  \cite[p.42, Eq.(1)]{AB36} (see also \cite{AB17}) with $p$ numerator parameters $\alpha_{1},  ... , \alpha_{p}$ and $q$ denominator parameters $\beta_{1}, ..., \beta_{q}$ defined by\\
\begin{equation}\label{h1}
{}_{p}F_{q}\left(\ \begin{array}{lll}\alpha_{1},...,\alpha_{p}~;~\\\beta_{1}, ..., \beta_{q}~;~\end{array} z\right)
=\sum_{n=0}^{\infty}\frac{(\alpha_{1})_{n}  ...  (\alpha_{p})_{n}}{(\beta_{1})_{n} ... (\beta_{q})_{n}}\frac{z^{n}}{n!}~,
\end{equation}
where $\alpha_{j}\in{\bf C}~(1\leq j\leq p)$ and $\beta_{j}\in{\bf C}\setminus {\bf Z}_{0}^{-}~(1\leq j\leq q)$ and $p,~q\in{\bf N}_{0}$.
The series ${}_{p}F_{q}(\cdot)$ is convergent for $|z|<\infty$ if $p\leq q$, and for $|z|<1$ if $p=q+1$.
Furthermore, if we set
\[
\omega=\sum_{j=1}^{q}\beta_{j}-\sum_{j=1}^{p}\alpha_{j},
\]
it is known that the ${}_{p}F_{q}$ series, with $p=q+1$, is absolutely convergent on the unit circle $|z|=1$ if $\Re(\omega)>0,$ and conditionally convergent on the unit circle $|z|=1, z\neq1$, if $-1<\Re(\omega)\leq0$.
%iii. divergent for $|z|=1$ if $Re(\omega)\leq-1$\\

The modified Bessel function of the second kind $K_{\nu}(z)$ of order $\nu$ (also known as the Macdonald function) is defined by (see \cite[p.~251]{AB26},  \cite{AB32})
\begin{eqnarray}\label{EB4}
K_{\nu}(z)=\sqrt{\pi}(2z)^{\nu} e^{-z}\,U(\nu+\fs, 2\nu+1, 2z),
\end{eqnarray}
where $U(a,b,z)$ is the confluent hypergeometric function \cite[p.~322]{AB26}. 
%Bateman's k-function \cite[p.40, ch.1, eq.(30)]{AB36}, is given by
%\begin{eqnarray}\label{ES18}
%k_{v}(z)
%=\frac{1}{\Gamma(1+\frac{v}{2})}W_{\frac{v}{2},\pm\frac{1}{2}}(2z),\\
%=\frac{e^{-z}}{\Gamma(1+\frac{v}{2})}U(-\frac{v}{2},0, 2z).
%\end{eqnarray}
The Meijer $G$-function is defined by means of the Mellin-Barnes contour integral \cite[(16.17.1)]{AB26}
% When $k=1,2,...,n$ and $\ell=1,2,...,m$, and $\alpha_{k}-\beta_{\ell}\neq$ positive integer then
\[
G_{{p} ,{q}}^{{m},{n}} \left( z~\bigg{|} \begin{array}{lll} \alpha_{1},...,\alpha_{n};\alpha_{n+1},...,\alpha_{p} \\\beta_{1},...,\beta_{m};\beta_{m+1},...,\beta_{q} \end{array} \right)\hspace{6cm}\]
\begin{equation}\label{RCG32}
\hspace{3cm}=\frac{1}{2\pi i}\int_{-i\infty}^{+i\infty}
\frac{\displaystyle\prod_{j=1}^{m}\Gamma(\beta_{j}-\zeta)\prod_{j=1}^{n}\Gamma(1-\alpha_{j}+\zeta)}
{\displaystyle\prod_{j=m+1}^{q}\Gamma(1-\beta_{j}+\zeta)\prod_{j=n+1}^{p}\Gamma(\alpha_{j}-\zeta)}z^{\zeta}d\zeta,
\end{equation}
where $z\neq0$, and $m,n,p,q$ are non-negative integers such that $ 1\leq m\leq q$~, $0\leq n\leq p$ and $p\leq q$.
The integral (\ref{RCG32}) converges in the sector $|\arg\,z|<\pi\kappa$, where $\kappa=m+n-\frac{1}{2}(p+q)$ and it is supposed that $\kappa>0$.
%If
%\begin{equation}\label{RCG33}
%\kappa=m+n-\frac{1}{2}(p+q).
%\end{equation}
%then the integral (\ref{RCG32}) converges in the sector $|\arg\,z|<\pi\kappa$ when $\kappa > 0$. 
%(ii). If $|\arg(z)|=\pi\kappa$, and $\kappa\geq0$, then the integral (\ref{RCG32}) converges absolutely.\\

The $G$-function is important in applied mathematics and formulas developed for the 
$G$-function become  master or key formulas from which a very large number of relations can be deduced for Bessel functions, their combinations and many other related functions. Thus the following list of some particular cases of  Meijer's $G$-function associated with the Bessel function $K_{\nu}(z)$ has been obtained mainly from several papers by C. S. Meijer  (see also \cite[pp.~219-220, Eq.(47-50)]{AB17}, \cite[p.48, Eq.(12)]{AB36})
\begin{equation}\label{ES19}
K_{\nu}(z)=\sqrt{\pi}~e^{z}~G_{{1} ,{2}}^{{2},{0}} \left( 2z \bigg{|} \begin{array}{lll} \frac{1}{2}\\ \nu~,~-\nu \end{array} \right),
\end{equation}
\begin{equation}\label{SE2}
=\frac{\cos (\pi\nu)}{\sqrt{\pi}}~e^{-z}~G_{{1} ,{2}}^{{2},{1}} \left( 2z \bigg{|} \begin{array}{lll} \frac{1}{2}\\ \nu~,~-\nu \end{array} \right),
\end{equation}
\begin{equation}\label{SE3}
=z^{-\mu}2^{\mu-1}G_{{0} ,{2}}^{{2},{0}} \left( \frac{z^{2}}{4} \bigg{|} \begin{array}{lll}  \\ \frac{\mu+\nu}{2}~,~\frac{\mu-\nu}{2} \end{array} \right),
\end{equation}
%\begin{equation}\label{SE4}
%=(2z)^{-\mu}e^{z}~G_{{1} ,{2}}^{{2},{1}} \left( 2z \bigg{|} \begin{array}{lll} \mu+\frac{1}{2} \\ \mu+\nu,~\mu-\nu \end{array} \right),
%\end{equation}
\begin{equation}\label{SE5}
=\cos(\pi\nu)\frac{(2z)^{-\mu}e^{z}}{\sqrt{\pi}}~ G_{{1} ,{2}}^{{2},{1}} \left( 2z \bigg{|} \begin{array}{lll} \mu+\frac{1}{2} \\ \mu+\nu,~\mu-\nu \end{array} \right),
\end{equation}
\begin{equation}\label{SE6}
=\frac{z^{-\mu}4^{\mu-1}}{\pi}G_{{0} ,{4}}^{{4},{0}} \left( \frac{z^{4}}{256} \bigg{|} \begin{array}{lll}  \\ \frac{\mu+\nu}{4}~,\frac{2+\mu+\nu}{4},~\frac{\mu-\nu}{4},\frac{2+\mu-\nu}{4} \end{array} \right),
\end{equation}
where $\mu$ is a free parameter and in all these expressions we have $z\neq 0$.
%In view of the above equations (\ref{ES19}) to (\ref{SE6}) where $z\neq0$.\\
%Also, other two results of particular cases of the Meijer's G-function associated with the Bessel function of second kind is given by [ref H.M.Sri]
In 1997, Chaudhry {\it et al.} \cite[Eq.(1.7)]{AB8} gave a $p$-extension of the Beta function $B(x, y)$ given by \[B(x,y; p)= \int_0^1 t^{x-1}(1-t)^{y-1}\,\exp \left[\frac{-p}{t(1-t)}\right]dt,\qquad (\Re (p)>0)\]
and they proved that this extension has connections with the Macdonald, error and Whittaker functions. Also,  Chaudhry {\it et al.} \cite{AB10} extended the Gaussian hypergeometric series ${}_{2}F_{1}(\cdot)$ and its integral representations. 
Recently, Parmar {\it et al.} \cite{AB30}, have given a further extension of the extended Beta function $B(x,y;p)$ by adding one more parameter $\nu$, which we denote and define by
\begin{equation}\label{ES20}
 B_{p,\nu}(x,y)=\sqrt{\frac{2p}{\pi}}\displaystyle{\int_{0}^{1}t^{x-\frac{3}{2}}(1-t)^{y-\frac{3}{2}}}
 K_{\nu+\frac{1}{2}}\left(\frac{p}{t(1-t)}\right)dt,
 \end{equation}
where $\Re(p)>0$, $\nu\geq0$ and $K_{\nu+\frac{1}{2}}(\cdot)$ is the modified Bessel function of order $\nu+\frac{1}{2}$. When $\nu=0$, (\ref{ES20}) reduces to $B(x,y;p)$, since $K_\frac{1}{2}(z)=(\pi/(2z))^\frac{1}{2} e^{-z}$. A different generalization of the Beta function has been given in \cite{AB28}.

%Mellin transforms of an integral function $f(x,y)$ over the interval $(0,\infty)$ with indices $r$ and $s$ is defined by
%\begin{equation}\label{ES04}
%{\cal M}\{f(x,y)\}(r,s)=\int_{0}^{\infty}\int_{0}^{\infty}x^{r-1}y^{s-1}f(x,y)dxdy
%\end{equation}
%The integral formula over the interval $(0,\infty)$ \cite[eq.(10.43.19)]{AB26}, is given by
%\begin{equation}\label{ES004}
%\int_{0}^{\infty}w^{s-\frac{1}{2}}K_{\alpha+\frac{1}{2}}(w)dw=2^{s-\frac{3}{2}}\Gamma\left(\frac{(s-\alpha)}{2}\right)
%\Gamma\left(\frac{(s+\alpha+1)}{2}\right).
%\end{equation}
%where $|\Re(\alpha+\frac{1}{2})|<\Re(s),~\Re(s-\alpha)>0,~\Re(s+\alpha)>-1$.\\
%Mellin transform representation of the new generalized beta function is given by
%\begin{equation}\label{ES400}
%\int_{0}^{\infty}p^{s-1}B_{p}^{(\alpha,\beta)}(x,y)dp=B(s+x,y+s)\Gamma^{(\alpha,\beta)}(s),
%\end{equation}
%where $\Re(s)>0; \Re(x + s)>0, \Re(y + s)>0; \Re(p)>0; \Re(\alpha)>0, \Re(\beta)>0$.\\
%It is obvious that $B_{p}^{(\alpha,\beta)}(x,y)=B_{0}^{(\alpha,\alpha)}(x,y)=B(x,y)$ and $ \Gamma^{(\alpha,\beta)}(s)=\Gamma^{(\alpha,\alpha)}(s)=\Gamma(s)$.\\
%In view of the above eq.(\ref{ES400}), we have
%\begin{equation}\label{ES401}
%\int_{0}^{\infty}p^{s-1}B(x,y)dp=B(s+x,y+s)\Gamma(s)
%\end{equation}
%where $ \Re(s)>0; \Re(x + s)>0, \Re(y + s)>0.$\\

%*****************************************

Motivated by some of the above-mentioned extensions of special functions, many authors have studied integral representations of $F_{1}(\cdot)$ functions.  Our aim in this paper is to introduce a $(p,v)$-extension of the Appell hypergeometric function in (\ref{ES9}), which we denote by $F_{1,p,\nu}(\cdot)$, based on the extended Beta function in (\ref{ES20}), and to systematically investigate some properties of this extended function.
We consider the Mellin transform, a differential formula, recursion formulas and a bounded inequality satisfied by this function. Also, we obtain some integral representations for $F_{1,p,v}$ containing Meijer's $G$-function. 
%and other contains Whittaker and confluent hypergeometric function. 

The plan of this paper as follows. The extended Appell function $F_{1,p,\nu}(\cdot)$  and its integral representation are defined in Section 2. Some new integral representations for $F_{1,p,v}(\cdot)$ involving the Meijer $G$-function are given.
The main properties $F_{1,p,v}(\cdot)$, namely its  Mellin transform, a differential formula, recurrence relation and a bounded inequality are established in Sections 3--6. Some concluding remarks are made in Section 7.

\vspace{0.6cm}

\begin{center}
{\bf 2. \  The $(p,\nu)$-extended Appell function $F_{1,p,\nu}(\cdot)$}
\end{center}
\setcounter{section}{2}
\setcounter{equation}{0}
\renewcommand{\theequation}{\arabic{section}.\arabic{equation}}
In \cite{AB27}, $\ddot{O}$zarslan {\it et al.} gave an extension of Appell's hypergeometric function  $F_{1}(\cdot)$ together with its integral representation. Here we consider the following $(p,v)$-extension of the Appell  hypergeometric function, which we denote by $F_{1,p,\nu}(\cdot)$, based on the extended beta function $B_{p,\nu}(x,y)$ defined in (\ref{ES20}). This is given by 
\begin{equation}\label{ES23}
F_{1,p,\nu}(b_{1},b_{2},b_{3};c_{1};x,y)=\sum_{n,m=0}^{\infty}\frac{(b_{2})_{m}(b_{3})_{n}~B_{p,\nu}(b_{1}+m+n,c_{1}-b_{1})}{~B(b_{1},c_{1}-b_{1})}\frac{x^{m}}{m!}\frac{y^{n}}{n!},
\end{equation}
where $|x|<1,~~|y|<1$ and $b_{1},b_{2},b_{3}\in{\bf C}$ and $c_{1}\in{\bf C}\backslash {\bf Z}_{0}^{-}$.
This definition clearly reduces to the original function when $\nu=0$.

An integral representation for the function $F_{1,p,\nu}(\cdot)$ is given by
\[F_{1,p,\nu}(b_{1},b_{2},b_{3};c_{1};x,y)=\frac{\Gamma(c_{1})}{\Gamma(b_{1})\Gamma(c_{1}-b_{1})}\times\\\]
\bee\label{ES24}
\times
\sqrt{\frac{2p}{\pi}}\int_{0}^{1}t^{b_{1}-\frac{3}{2}}(1-t)^{c_{1}-b_{1}-\frac{3}{2}}(1-xt)^{-b_{3}}(1-yt)^{-b_{2}}
K_{\nu+\frac{1}{2}}\left(\frac{p}{t(1-t)}\right)dt,
\ee
where $\Re(p)>0$, $\nu\geq0$, $|\arg(1-x)|<\pi$ and $|\arg(1-y)|<\pi$ and we impose the condition $\Re(c_1)>\Re(b_1)>0$ for the multiplicative factor $1/B(b_1,c_1-b_1)$ to be finite..
That this representation yields (\ref{ES23}) can be shown by binomially expanding the factors $(1-xt)^{-b_3}$ and $(1-yt)^{-b_2}$ when $|x|$, $|y|<1$, reversing the order of summation and integration and evaluating the resulting integral by (\ref{ES20}).
%The domain of convergence for the extended Appell hypergeometric function $F_{1,p,\nu}(\cdot)$ is given in \cite[p.4]{AB}.

%
\newtheorem{theorem}{Theorem}
\begin{theorem}$\!\!\!.$\ Each of the following integral representations of $F_{1,p,\nu}(\cdot)$ associated with Meijer's $G$-function holds for $p\in{\bf C}_{>}$.
\[F_{1,p,v}(b_{1},b_{2},b_{3};c_{1};x,y)=\frac{\Gamma(c_{1})\sqrt{2p}}{\Gamma(b_{1})\Gamma(c_{1}-b_{1})}
\int_{0}^{1}t^{b_{1}-\frac{3}{2}}(1-t)^{c_{1}-b_{1}-\frac{3}{2}}(1-xt)^{-b_{3}}(1-yt)^{-b_{2}}\times\]
\bee\label{ES65}
 \times
e^{\frac{p}{t(1-t)}} G_{{1} ,{2}}^{{2},{0}} \left( \frac{2p}{t(1-t)} \bigg{|} \begin{array}{lll}\frac{1}{2}
\\ \nu+\frac{1}{2}~,-\nu-\frac{1}{2} \end{array} \right)dt,
\ee
%where $\left(\Re(c_{1})>\Re(b_{1})>0, |arg(1-x)|<\pi, |arg(1-y)|<\pi\right)$,
\[=\frac{\Gamma(c_{1})\sqrt{2p}}{\Gamma(b_{1})\Gamma(c_{1}-b_{1})}\frac{\cos\pi(\nu\!+\!\fs)}{\pi}
\int_{0}^{1}t^{b_{1}-\frac{3}{2}}(1-t)^{c_{1}-b_{1}-\frac{3}{2}}(1-xt)^{-b_{3}}\times\]
\bee\label{ES66}
 \times(1-yt)^{-b_{2}}
e^{\frac{-p}{t(1-t)}} G_{{1} ,{2}}^{{2},{1}} \left( \frac{2p}{t(1-t)} \bigg{|} \begin{array}{lll}\frac{1}{2}
\\ \nu+\frac{1}{2}~,-\nu-\frac{1}{2} \end{array} \right)dt,
\ee
%where $\left(\Re(c_{1})>\Re(b_{1})>0, |arg(1-x)|<\pi, |arg(1-y)|<\pi\right)$,
\[=\frac{\Gamma(c_{1})2^{\mu-\frac{1}{2}}p^{-\mu+\frac{1}{2}}}{\Gamma(b_{1})\Gamma(c_{1}-b_{1})}
\int_{0}^{1}t^{b_{1}+\mu-\frac{3}{2}}(1-t)^{c_{1}-b_{1}+\mu-\frac{3}{2}}(1-xt)^{-b_{3}}\times\]
\bee\label{ES67}
 \times(1-yt)^{-b_{2}}
G_{{0} ,{2}}^{{2},{0}} \left( \frac{p^{2}}{4t^{2}(1-t)^{2}} \bigg{|} \begin{array}{lll}  \\ \frac{2\mu+2\nu+1}{4}~,~\frac{2\mu-2\nu-1}{4} \end{array} \right)dt,
\ee
%where $\left(\Re(c_{1})>\Re(b_{1})>0, |arg(1-x)|<\pi, |arg(1-y)|<\pi\right)$,
\[=\frac{\Gamma(c_{1})(2p)^{-\mu+\frac{1}{2}}}{\Gamma(b_{1})\Gamma(c_{1}-b_{1})}\frac{\cos \pi(\nu\!+\!\fs)}{\pi}
\int_{0}^{1}t^{b_{1}+\mu-\frac{3}{2}}(1-t)^{c_{1}-b_{1}+\mu-\frac{3}{2}}(1-xt)^{-b_{3}}\times\]
\bee\label{ES68}
 \times(1-yt)^{-b_{2}}
e^{\frac{-p}{t(1-t)}}G_{{1} ,{2}}^{{2},{1}} \left( \frac{2p}{t(1-t)} \bigg{|} \begin{array}{lll} \frac{2\mu+1}{2} \\ \frac{2\mu+2\nu+1}{2}~,~\frac{2\mu-2\nu-1}{2} \end{array} \right)dt,
\ee
%where $\left(\Re(c_{1})>\Re(b_{1})>0, |arg(1-x)|<\pi, |arg(1-y)|<\pi\right)$,
\[=\frac{\Gamma(c_{1})p^{-\mu+\frac{1}{2}}2^{2\mu-\frac{3}{2}}}{\Gamma(b_{1})\Gamma(c_{1}-b_{1})\pi^{\frac{3}{2}}}
\int_{0}^{1}t^{b_{1}+\mu-\frac{3}{2}}(1-t)^{c_{1}-b_{1}+\mu-\frac{3}{2}}(1-xt)^{-b_{3}}\times\]
\bee\label{ES69}
 \times(1-yt)^{-b_{2}}
G_{{0} ,{4}}^{{4},{0}} \left( \frac{p^{4}}{(4t)^{4}(1-t)^{4}} \bigg{|} \begin{array}{lll}  \\ \frac{2\mu+2\nu+1}{8}~,~\frac{2\mu+2\nu+5}{8},\frac{2\mu-2\nu-1}{8}~,~\frac{2\mu-2\nu+3}{8} \end{array} \right)dt,
\ee
where $\Re(c_{1})>\Re(b_{1})>0, |\arg(1-x)|<\pi, |\arg(1-y)|<\pi$ and $\mu$ is a free parameter.
\end{theorem}
\noindent
\textbf{Proof}: The above integral representations (\ref{ES65})--(\ref{ES69}) are obtained by using (\ref{ES19})--(\ref{SE6}) in the expression of the extended Appell function in (\ref{ES24}).
Similarly, other new integral representations of $F_{1,p,v}(\cdot)$ associated with the confluent hypergeometric function can be obtained using (\ref{EB4}) in (\ref{ES24}).

The following transformation formula can be derived from the integral representation (\ref{ES24}) for $F_{1,p,\nu}(\cdot)$.
\begin{theorem}$\!\!\!.$\ The following transformation formula holds:
\bee\label{e21}
F_{1,p,\nu}(b_1,b_2,b_3;c_1;x,y)=(1-x)^{-b_3}(1-y)^{-b_2} F_{1,p,\nu}\left(c_1-b_1,b_2,b_3;c_1;\frac{x}{x-1},\frac{y}{y-1}\right).
\ee
\end{theorem}
\noindent
\textbf{Proof}: Put $t=1-\zeta$ in (\ref{ES24}) to obtain
\[F_{1,p,\nu}(b_1,b_2,b_3;c_1;x,y)=\frac{\Gamma(c_1)(1-x)^{-b_3}(1-y)^{-b_2}}{\Gamma(b_1) \Gamma(c_1-b_1)}\times\]
\[\sqrt{\frac{2p}{\pi}} \int_0^1 \zeta^{c_1-b_1-\frac{3}{2}}(1-\zeta)^{b_1-\frac{3}{2}}\left(1-\frac{x}{x-1}\zeta\right)^{-b_3}\left(1-\frac{y}{y-1}\zeta\right)^{-b_2} K_{\nu+\frac{1}{2}}\left(\frac{p}{\zeta(1-\zeta)}\right)d\zeta.\]
Identification of the above integral as a $F_{1,p,\nu}(\cdot)$ function then yields (\ref{e21}).

\vspace{0.6cm}

\begin{center}
{\bf 3. \  The Mellin transform of $F_{1,p,\nu}(\cdot)$}
\end{center}
\setcounter{section}{3}
\setcounter{equation}{0}
\renewcommand{\theequation}{\arabic{section}.\arabic{equation}}
The Mellin transform of a locally integrable function $f(x)$ on $(0,\infty)$ is defined by
\[{\cal M}\{f(x)\}(s)=\int_0^\infty x^{s-1}f(x)\,dx\]
when the integral converges.\\
%We have the following theorem:
%
\begin{theorem}$\!\!\!.$\ The following Mellin transform of the extended Appell hypergeometric function $ F_{1,p,v}(\cdot)$ holds true:
\begin{eqnarray}\label{ES70}
{\cal M}\left\{F_{1,p,v}(b_{1},b_{2},b_{3};c_{1};x,y)\right\}(s)=\int_{0}^{\infty}p^{s-1}F_{1,p,v}(b_{1},b_{2},b_{3};c_{1};x,y)\,dp\nonumber\\
=\frac{2^{s-1}}{\sqrt{\pi}}\Gamma\left(\frac{s-v}{2}\right)\Gamma\left(\frac{s+v+1}{2}\right)~F_{1}\left(b_{1}+s,b_{2},b_{3};c_{1}+s;x,y\right),
\end{eqnarray}
where $\Re(s-v)>0$, $\Re(s+v)>-1$, $\Re(s)>0$ and $c_{1}+s\in{\bf C}/{\bf Z}_{0}^{-}$.
\end{theorem}
\noindent
\textbf{Proof}: Substituting the extended Appell function (\ref{ES23}) into the left-hand side of (\ref{ES70}) and changing the order of integration (by the uniform convergence of the integral), we obtain
\[{\cal M}\left\{F_{1,v}(b_{1},b_{2},b_{3};c_{1};x,y;p)\right\}(s)
= \frac{\Gamma(c_{1})}{\Gamma(b_{1})\Gamma(c_{1}-b_{1})}
\sqrt{\frac{2}{\pi}}\int_{0}^{1}t^{b_{1}-\frac{3}{2}}(1-t)^{c_{1}-b_{1}-\frac{3}{2}}(1-xt)^{-b_{3}}\times\]
\bee\label{ES71}
\times(1-yt)^{-b_{2}}
\left\{\int_{0}^{\infty}p^{s-\frac{1}{2}}K_{v+\frac{1}{2}}\left(\frac{p}{t(1-t)}\right)dp\right\}\,dt.
\ee
Application of the result \cite[(10.43.19)]{AB26}
\[\int_{0}^{\infty}w^{s-\frac{1}{2}}K_{\alpha+\frac{1}{2}}(w)dw=2^{s-\frac{3}{2}}\Gamma\left(\frac{s-\alpha}{2}\right)
\Gamma\left(\frac{s+\alpha+1}{2}\right)\qquad (|\Re (\alpha)|<\Re (s))\]
followed by the substitution $w=p/\{t(1-t)\}$ in (\ref{ES71}) then yields
\[\varphi(s)\equiv{\cal M}\left\{F_{1,p,v}(b_{1},b_{2},b_{3};c_{1};x,y)\right\}(s)
= \frac{2^{s-1}}{\sqrt{\pi}}\Gamma\left(\frac{s-v}{2}\right)\Gamma\left(\frac{s+v+1}{2}\right)\times\]
\bee\label{ES72}
 \times \frac{\Gamma(c_{1})}{\Gamma(b_{1})\Gamma(c_{1}-b_{1})}
\int_{0}^{1}t^{b_{1}+s-1}(1-t)^{c_{1}+s-b_{1}-1}(1-xt)^{-b_{3}}(1-yt)^{-b_{2}}\,dt.
\ee
Finally, using the definition of the Appell function $F_1(\cdot)$ in (\ref{ES10}), we obtain the right-hand side of (\ref{ES70}).
\bigskip

\noindent
\textbf{Corollary}: \ The following inverse Mellin formula for $ F_{1,p,v}(\cdot)$ holds:
\[
F_{1,p,v}(b_{1},b_{2},b_{3};c_{1};x,y)={\cal M}^{-1}\left\{\varphi(s)\right\}\]
\bee
=\frac{1}{4\pi i\sqrt{\pi}}\int_{c-i\infty}^{c+i\infty}\left(\frac{2}{p}\right)^{s}\Gamma\left(\frac{s-v}{2}\right)\Gamma\left(\frac{s+v+1}{2}\right)
F_{1}\left(b_{1}+s,b_{2},b_{3};c_{1}+s;x,y\right)\,ds,
\ee
where $c>\nu$.
\vspace{0.6cm}

\begin{center}
{\bf 4. \  A differentiation formula for $F_{1,p,\nu}(\cdot)$}
\end{center}
\setcounter{section}{4}
\setcounter{equation}{0}
\renewcommand{\theequation}{\arabic{section}.\arabic{equation}}
\begin{theorem}$\!\!\!.$\ The following differentiation formula for $F_{1,p,v}(\cdot)$ holds:
\[\frac{\partial^{M+N}}{\partial x^{M}\partial y^{N}}F_{1,p,\nu}(b_{1},b_{2},b_{3};c_{1};x,y)\hspace{6cm}\]
\bee\label{ES75}
=\frac{(b_{1})_{M+N}(b_{2})_M(b_{3})_N}{(c_{1})_{M+N}}F_{1,p,\nu}(b_{1}+M+N,b_{2}+M,b_{3}+N;c_{1}+M+N;x,y),
\ee
where $M,N\in{\bf N}_{0}$.
\end{theorem}
\noindent
\textbf{Proof}:\ If we differentiate the series for $F_{1,p,\nu}(b_{1},b_{2},b_{3};c_{1};x,y)$ in (\ref{ES23}) with respect to $x$ we obtain
\[\frac{\partial}{\partial x} F_{1,p,\nu}(b_{1},b_{2},b_{3};c_{1};x,y)
=\sum_{m=1}^\infty \sum_{n=0}^\infty \frac{(b_{2})_{m}(b_{3})_{n}~B_{p,\nu}(b_{1}+m+n,c_{1}-b_{1})}{~B(b_{1},c_{1}-b_{1})}\frac{x^{m-1}}{(m-1)!}\frac{y^{n}}{n!}.\]
Making use of the fact that 
\bee\label{e40}
B(b_1,c_1-b_1)=\frac{c_1}{b_1}\,B(b_1+1, c_1-b_1)
\ee
and $(b_2)_{m+1}=b_2(b_2+1)_m$, we have upon setting $m \to m+1$
\begin{eqnarray}
\frac{\partial}{\partial x} F_{1,p,\nu}(b_{1},b_{2},b_{3};c_{1};x,y)\!\!&=&\!\!\frac{b_1b_2}{c_1} \sum_{m,n=0}^\infty 
\frac{(b_{2}+1)_{m}(b_{3})_{n}~B_{p,\nu}(b_{1}+1+m+n,c_{1}-b_{1})}{~B(b_{1}+1,c_{1}-b_{1})}\frac{x^{m}}{m!}\frac{y^{n}}{n!}\nonumber\\
&=&\!\!\frac{b_1b_2}{c_1}\,F_{1,p,\nu}(b_1+1,b_2+1,b_3;c_1+1;x,y).\label{e41}
\end{eqnarray}
Repeated application of (\ref{e41}) then yields for $M=1, 2, \ldots $
\[\frac{\partial^M}{\partial x^M} F_{1,p,\nu}(b_{1},b_{2},b_{3};c_{1};x,y)=
\frac{(b_1)_M (b_2)_M}{(c_1)_M}\,F_{1,p, \nu}(b_1+M,b_2+M,b_3;c_1+M;x,y).\]

A similar reasoning shows that
\[\frac{\partial^{M+1}}{\partial x^M \partial y}F_{1,p,\nu}(b_{1},b_{2},b_{3};c_{1};x,y)\hspace{5cm}\]
\[=\frac{(b_1)_M (b_2)_M}{(c_1)_M} \sum_{m=0}^\infty \sum_{n=1}^\infty \frac{(b_2+M)_m (b_3)_n B_{p, \nu}(b_1+M+m+n,c_1-b_1)}{B(b_1+M, c_1-b_1)}\,\frac{x^m}{m!} \frac{y^{n-1}}{(n-1)!}\]
\bee\label{e42}
=\frac{(b_1)_{M+1} (b_2)_M b_3}{(c_1)_{M+1}}\,F_{1,p, \nu}(b_1+M+1,b_2+M,b_3+1;c_1+M+1;x,y).
\ee
Repeated differentiation of (\ref{e42}) $N$ times with respect to $y$ then readily produces the result stated in (\ref{ES75}).
The result (\ref{ES75}) has been derived assuming that $|x|<1$, $|y|<1$ but can be extended to all values of $x$ and $y$ satisfying $|\arg (1-x)|<\pi$, $|\arg (1-y)|<\pi$ by appeal to analytic continuation.
\vspace{0.6cm}

\begin{center}
{\bf 5. \  An upper bound for $F_{1,p,\nu}(\cdot)$}
\end{center}
\setcounter{section}{5}
\setcounter{equation}{0}
\renewcommand{\theequation}{\arabic{section}.\arabic{equation}}
\begin{theorem}$\!\!\!.$\ Let the parameters $b_1$, $b_2$, $b_3$, $c_1$ and the variables $x$, $y$ be real. Then the following bounded inequality for $ F_{1,p,\nu}(\cdot)$ holds:
\[|F_{1,p,\nu}(b_{1},b_{2},b_{3};c_{1};x,y)|\hspace{8cm}\]
\bee\label{e51}
<\frac{2^\nu |p|^{\nu+1}}{\sqrt{\pi} (\Re (p))^{2\nu+1}}\,\frac{\Gamma(\nu+\fs)B(b_1+\nu,c_1-b_1+\nu)}{B(b_1,c_1-b_1)}\,F_1(b_1+\nu,b_2,b_3;c_1+2\nu;x,y),
\ee
where $\Re(p)>0$.
\end{theorem}

The integral representation of the extension $ F_{1,p,\nu}(\cdot)$ in (\ref{ES24}) is associated with the modified Bessel function of the second kind, for which we have the following expression \cite[(10.32.8)]{AB26}
\[ K_{\nu+\frac{1}{2}}(z)=\frac{\sqrt{\pi}\left(\frac{1}{2}z\right)^{\nu+\frac{1}{2}}}{\Gamma(\nu+1)}\int_{1}^{\infty}e^{-zt}(t^{2}-1)^{\nu}dt,\qquad(\nu>-1,\ \Re (z)>0).
\]
In our problem we have $\nu>0$, $\Re(z)>0$. Further, we let $x=\Re(z)$, so that
\[| K_{\nu+\frac{1}{2}}(z)|\leq\frac{\sqrt{\pi}\left(\frac{1}{2}|z|\right)^{\nu+\frac{1}{2}}}{\Gamma(\nu+1)}\left|\int_{1}^{\infty}e^{-zt}(t^{2}-1)^{\nu}dt\right|<\frac{\sqrt{\pi}\left(\frac{1}{2}|z|\right)^{\nu+\frac{1}{2}}}{\Gamma(\nu+1)} \int_0^1 t^{2\nu}e^{-xt} dt \]
\begin{equation}\label{ES87}
=\frac{\sqrt{\pi}\left(\frac{1}{2}|z|\right)^{\nu+\frac{1}{2}}}{\Gamma(\nu+1)}\frac{\Gamma(2\nu+1,x)}{x^{2\nu+1}},
\end{equation}
where $\Gamma(a,z)$ is the upper incomplete gamma function \cite[(8.2.2)]{AB26}. Although this bound is numerically found to be quite sharp when $z $ is real, it involves the incomplete gamma function which would make the integral for $F_{1,p,\nu}(b_{1},b_{2},b_{3};c_{1};x,y)$ difficult to bound. We can simplify (\ref{ES87}) by making use of the simple inequality $\Gamma(2\nu+1,x)<\Gamma(2\nu+1)$ to find
\begin{equation}\label{ES88}
| K_{\nu+\frac{1}{2}}(z)|<\frac{\sqrt{\pi}\left(\frac{1}{2}|z|\right)^{\nu+\frac{1}{2}}}{\Gamma(\nu+1)}\frac{\Gamma(2\nu+1)}{x^{2\nu+1}}
=\frac{1}{2}\left(\frac{2|z|}{x^{2}}\right)^{\nu+\frac{1}{2}}\Gamma(\nu+\fs),
\end{equation}
upon use of the duplication formula for the gamma function. The bound (\ref{ES88}) is less sharp than (\ref{ES87}) but has the advantage of being easier to handle in the integral for $F_{1,p,\nu}(b_{1},b_{2},b_{3};c_{1};x,y)$. 
\\
\\
\noindent
\textbf{Proof}: Setting $z=p/(t(1-t))$, where $t\in(0,1)$ and $\Re(p)>0$, in (\ref{ES88}) we obtain
\[
\left| K_{\nu+\frac{1}{2}}\left(\frac{p}{t(1-t)}\right)\right|<\frac{1}{2}\left(\frac{2|p|t(1-t)}{(\Re (p))^{2}}\right)^{\nu+\frac{1}{2}}\Gamma(\nu+\fs).
\]
For ease of presentation we shall assume that the parameters $b_1$, $b_2$, $b_3$ and $c_1$ are real; the extension to complex parameters is straightforward. In addition, we shall consider only real values of the variables $x$ and $y$. Then, from (\ref{ES24}),
\[|F_{1,p,\nu}(b_1 b_2, b_3;c_1;x,y)|\hspace{8cm}\]
\[\leq \frac{\sqrt{2|p|/\pi}}{B(b_1,c_1-b_1)}\int_0^1\left|t^{b_1-\frac{3}{2}}(1-t)^{c_1-b_1-\frac{3}{2}} (1-xt)^{-b_3} (1-yt)^{-b_2} \,K_{\nu+\frac{1}{2}}\left(\frac{p}{t(1-t)}\right) \right|\,dt\]
\[<\frac{2^\nu |p|^{\nu+1}}{\sqrt{\pi} (\Re (p))^{2\nu+1}}\,\frac{\Gamma(\nu+\fs)}{B(b_1,c_1-b_1)} \int_0^1
t^{b_1+\nu-1}(1-t)^{c_1-b_1+\nu-1} (1-xt)^{-b_3} (1-yt)^{-b_2}dt\]
\bee\label{e52}
<\frac{2^\nu |p|^{\nu+1}\Gamma(\nu+\fs)}{\sqrt{\pi} (\Re (p))^{2\nu+1}}\,\frac{B(b_1+\nu,c_1-b_1+\nu)}{B(b_1,c_1-b_1)}\,F_1(b_1+\nu,b_2,b_3;c_1+2\nu;x,y),
\ee
which is the result stated in (\ref{e51}).
 
If we have $x<0$, $y<0$ (resp. $x>0$, $y>0$) and suppose further that $b_2>0$, $b_3>0$ (resp. $b_2<0$, $b_3<0$) then we obtain the simpler bound
\[|F_{1,p,\nu}(b_1 b_2, b_3;c_1;x,y)|<\frac{2^\nu |p|^{\nu+1}\Gamma(\nu+\fs)}{\sqrt{\pi} (\Re (p))^{2\nu+1}}\,\frac{B(b_1+\nu,c_1-b_1+\nu)}{B(b_1,c_1-b_1)}.\]
 \vspace{0.6cm}

\begin{center}
{\bf 6. \  Recursion formulas for $F_{1,p,\nu}(\cdot)$}
\end{center}
\setcounter{section}{6}
\setcounter{equation}{0}
\renewcommand{\theequation}{\arabic{section}.\arabic{equation}}
In view of the recursion formulas for the Appell function $F_1(\cdot)$ (see \cite{AB00} and \cite{AB42}) we give the following recursion formulas for the extended Appell function $F_{1,p,\nu}(\cdot)$.
\begin{theorem}$\!\!\!.$\ The following recursion formulas for the extended Appell function with respect to the numerator parameters $b_2$ and $b_3$ hold:
\[F_{1,p,\nu}(b_1,b_2+n,b_3;c_1;x,y)\hspace{10cm}\]
\bee\label{e61}
=F_{1,p,\nu}(b_1,b_2,b_3;c_1;x,y)+\frac{b_1x}{c_1} \sum_{\ell=1}^n F_{1,p,\nu}(b_1+1,b_2+\ell,b_3;c_1+1;x,y)
\ee
and
\[F_{1,p,\nu}(b_1,b_2,b_3+n;c_1;x,y)\hspace{10cm}\]
\bee\label{e62}
=F_{1,p,\nu}(b_1,b_2,b_3;c_1;x,y)+\frac{b_1y}{c_1} \sum_{\ell=1}^n F_{1,p,\nu}(b_1+1,b_2,b_3+\ell;c_1+1;x,y)
\ee
for positive integer $n$.
\end{theorem}

\noindent
\textbf{Proof}: From (\ref{ES23}) and the result $(b_2+1)_m=(b_2)_m (1+m/b_2)$, we obtain
\[
F_{1,p,\nu}(b_1,b_2+1,b_3;c_1;x,y)=\sum_{m,n=0}^\infty \frac{(b_2+1)_m(b_3)_n B_{p,\nu}(b_1+m+n,c_1-b_1)}{B(b_1,c_1-b_1)}\,\frac{x^m}{m!}\,\frac{y^n}{n!}\]
\[=F_{1,p,\nu}(b_1,b_2,b_3;c_1;x,y)+\frac{x}{b_2}\sum_{m=1}^\infty\sum_{n=0}^\infty
\frac{(b_2)_m(b_3)_n B_{p,\nu}(b_1+m+n,c_1-b_1)}{B(b_1,c_1-b_1)}\,\frac{x^{m-1}}{(m\!-\!1)!}\,\frac{y^n}{n!}.\]
Setting $m\to m+1$ and using $(b_2)_{m+1}=b_2(b_2+1)_m$ together with (\ref{e40}), we find 
\[F_{1,p,\nu}(b_1,b_2+1,b_3;c_1;x,y)\hspace{10cm}\]
\[=F_{1,p,\nu}(b_1,b_2,b_3;c_1;x,y)+\frac{b_1x}{c_1}\sum_{m,n=0}^\infty \frac{(b_2+1)_m(b_3)_n B_{p,\nu}(b_1+1+m+n,c_1-b_1)}{B(b_1+1,c_1-b_1)}\,\frac{x^m}{m!}\,\frac{y^n}{n!}\]
\bee\label{e63}
=F_{1,p,\nu}(b_1,b_2,b_3;c_1;x,y)+\frac{b_1x}{c_1}\,F_{1,p,\nu}(b_1+1,b_2+1,b_3;c_1+1;x,y).
\ee

From (\ref{e63}) we obtain, upon putting $b_2\to b_2+1$,
\[F_{1,p,\nu}(b_1,b_2+2,b_3;c_1;x,y)\hspace{10cm}\]
\[=F_{1,p,\nu}(b_1,b_2+1,b_3;c_1;x,y)+\frac{b_1x}{c_1}\,F_{1,p,\nu}(b_1+1,b_2+2,b_3;c_1+1;x,y)\]
\[=F_{1,p,\nu}(b_1,b_2,b_3;c_1;x,y)+\frac{b_1x}{c_1}\sum_{\ell=1}^2 F_{1,p,\nu}(b_1+1,b_2+\ell,b_3;c_1+1;x,y).\]
Repeated application of the recursion (\ref{e63}) in this manner then immediately leads to the result stated in (\ref{e61}).
The proof of (\ref{e62}) is obtained in the same way by interchanging $b_2$ and $b_3$.

\vspace{0.6cm}

\begin{center}
{\bf 7. \  Concluding remarks}
\end{center}
\setcounter{section}{7}
\setcounter{equation}{0}
\renewcommand{\theequation}{\arabic{section}.\arabic{equation}}
%In the existing literature of special functions, integral representations of the special functions of applied mathematics and mathematical physics have been investigated. 
We have introduced the extended Appell hypergeometric function $ F_{1,p,\nu}(\cdot)$ given in (\ref{ES23})
by use of the extended Beta function defined in (\ref{ES20}).  Also, we have described some properties of this function, namely the Mellin transform, a differential formula, some recurrence relations and a bounded inequality. In addition, we have also obtained some new integral representations of the extended Appell hypergeometric function involving Meijer's $G$-function and indicated other possible representations in terms of the confluent hypergeometric function.

\end{document}